\documentclass{elsart}
\usepackage{amssymb}
\usepackage{amsmath}
\journal{Journal of Algebra}
\date{}

\newcommand\Qmax{Q^r_{\mathrm{max}}}
\newcommand\Qlmax{Q^l_{\mathrm{max}}}
\newcommand\Qlrmax{Q_{\mathrm{max}}}
\newcommand\Qtot{Q^r_{\mathrm{tot}}}
\newcommand\Qltot{Q^l_{\mathrm{tot}}}
\newcommand\Qlrtot{Q_{\mathrm{tot}}}
\newcommand\Qcl{Q^r_{\mathrm{cl}}}
\newcommand\Qlrcl{Q_{\mathrm{cl}}}
\newcommand\f{{\mathcal F}}
\newcommand\te{{\mathcal T}}
\newcommand\ef{{\mathfrak F}}
\def\dirlim{\mathop{\varinjlim}\limits}
\newcommand\ce{{\mathcal C}}
\newcommand\cl{\mathrm{cl}}
\newcommand\tor{\mathrm{Tor}}
\newcommand\homo{\mathrm{Hom}}

\begin{document}
\begin{frontmatter}

\title{A Simplification of Morita's Construction of Total Right
Rings of Quotients for a Class of Rings}

\author{Lia Va\v s}

\address{Department of Mathematics, Physics and Computer Science,
University of the Sciences in Philadelphia, 600 S. 43rd St.,
Philadelphia, PA 19104}

\ead{l.vas@usip.edu}

\begin{abstract} The total right ring of
quotients $Q_{\mathrm{tot}}^r(R),$ sometimes also called the
maximal flat epimorphic right ring of quotients or right flat
epimorphic hull, is usually obtained as a directed union of a
certain family of extension of the base ring $R$. In
\cite{Morita3}, $Q_{\mathrm{tot}}^r(R)$ is constructed in a
different way, by transfinite induction on ordinals. Starting with
the maximal right ring of quotients $Q_{\mathrm{max}}^r(R)$, its
subrings are constructed until $Q_{\mathrm{tot}}^r(R)$ is
obtained.

Here, we prove that Morita's construction of
$Q_{\mathrm{tot}}^r(R)$ can be simplified for rings satisfying
condition (C) that every subring of the maximal right ring of
quotients $Q^r_{\mathrm{max}}(R)$ containing $R$ is flat as a left
$R$-module. We illustrate the usefulness of this simplification by
considering the class of right semihereditary rings all of which
satisfy condition (C). We prove that the construction stops after
just one step and we obtain a simple description of
$Q^r_{\mathrm{tot}}(R)$ in this case. Lastly, we study conditions
that imply that Morita's construction ends in countably many
steps.
\end{abstract}


\end{frontmatter}

\section{Introduction}

There have been many attempts in ring theory to extend a given
ring $R$ to a ring in which some kind of generalized division is
possible. The classical right ring of quotients $\Qcl(R)$
unfortunately does not exist for every ring $R$. For many
important cases, the maximal right ring of quotients $\Qmax(R)$
always exists and has properties that bring it closer to being a
division ring. However, $\Qmax(R)$ may fail to have some
properties of $\Qcl(R)$ that we would prefer to keep.

Yet another attempt to find a reasonable right ring of quotients
was to consider the total right ring of quotients $\Qtot(R)$
sometimes also called the maximal flat epimorphic right ring of
quotients, right flat epimorphic hull or the maximal perfect right
localization. It can be defined for every ring and it is contained
in the maximal right ring of quotients. If the classical right
ring of quotients exists, the total right ring of quotients is
between the classical and the maximal right ring of quotients.
$\Qtot(R)$ is a generalization of the classical right ring of
quotients in the sense that every element $a\in\Qtot(R)$ has the
property
\[ar_i\in R\mbox{ and }\sum_{i=1}^n r_i a_i=1\mbox{ for some }n, a_i\in \Qtot(R)\mbox{ and }r_i\in R,\;i=1,\ldots,
n.\] Note that the above property implies that
\[a=a 1 = \sum_{i=1}^n  a r_i a_i = \sum_{i=1}^n s_i a_i\mbox{ where }s_i=ar_i\in R\mbox{ and }\sum_{i=1}^n r_i
a_i=1,\] which for $n=1,$ is the familiar property of the
classical right ring of quotients: every $a\in \Qcl(R)$ is of the
form $a=bt$ for some $b\in R$ and $t\in\Qcl(R)$ such that $t$ is
invertible in $\Qcl(R).$

Usually, the total right ring of quotients is constructed in the
following way. For any $R$, let us consider the family of all ring
extensions $S$ such that $S$ is flat as left $R$-module and that
the inclusion $R\subseteq S$ is an epimorphism in the category of
rings. This family is directed under inclusion. The directed union
of the elements of this family is the total right ring of
quotients $\Qtot(R).$ Several authors proved the existence of
$\Qtot(R)$ in a series of papers published in the late 1960s and
early 1970s: Findlay \cite{Findlay}, Knight \cite{Knight}, Lazard
\cite{Lazard}, Popescu and Spircu \cite{PopSpirc}. A good overview
of the subject is given in Stenstr\"om's book \cite{Stenstrom}.
Morita in \cite{Morita2} and \cite{Morita3} has a different
approach for defining $\Qtot(R).$ His idea is to start from the
maximal right ring of quotients $\Qmax(R)$ and to construct
$\Qtot(R)$ by transfinite induction on ordinals, "descending" from
$\Qmax(R)$ towards $R$ instead of "going upwards" starting from
$R$ using the directed family as in the classical construction.
This construction is described in the alternative proof of
Corollary 3.4 in \cite{Morita3}.

In this paper, we prove that Morita's construction of
$Q_{\mathrm{tot}}^r(R)$ can be simplified for rings that satisfy
the following condition
\begin{itemize}
\item[(C)] Every subring of the maximal right ring of quotients
$Q^r_{\mathrm{max}}(R)$ containing $R$ is flat as left $R$-module.
\end{itemize}

All rings constructed inductively in Morita's construction are
rings of right quotients of a certain torsion theory. The
simplification of the construction reduces to the simplification
of the description of this torsion theory. If the construction
ends after a finite number of steps, we obtain an explicit
description of $\Qtot(R)$.

A right semihereditary ring $R$ satisfies condition (C). We show
that the construction of $\Qtot(R)$ stops after at most one step
if $R$ is right semihereditary, producing the following
description of $\Qtot(R).$ An element $a$ of $\Qmax(R)$ is in
$\Qtot(R)$ if and only if
\[ar_i\in R\mbox{ and }\sum_{i=1}^n r_i a_i=1\mbox{ for some }n, a_i\in \Qmax(R)\mbox{ and }r_i\in R,\;i=1,\ldots, n.\]

In Section \ref{Section_general_quotients}, we review some basic
notions including torsion theories and right rings of quotients of
hereditary torsion theories. We also recall the definition and
basic properties of perfect right rings of quotients and the total
right ring of quotients. The exposition of rings of quotients
follows the one in \cite{Stenstrom}. This approach was first
introduced by Gabriel (see \cite{Gabriel}). In this section we
also present details of Morita's construction of $\Qtot(R).$

Section \ref{Section_construction} contains the construction of
$\Qtot(R)$ of a ring $R$ satisfying conditions (C). In Proposition
\ref{Morita=mine_whenC}, we prove that this construction and
Morita's coincide if $R$ satisfies condition (C).

In Section \ref{Section_semihereditary}, we turn our attention to
the class of right semihereditary rings and prove that the
construction ends after at most one step (Theorem
\ref{Qtot_for_semihereditary}). We illustrate the construction
with examples and survey the results on the condition that
Morita's construction ends already at the zeroth step.

In Section \ref{Section_construction_stops}, we study conditions
implying that the construction ends after countably many steps
(Proposition \ref{C_and_C'}).

We finish the paper by listing some interesting questions.

\section{Right Rings of Quotients}
\label{Section_general_quotients}

\subsection{General Right Rings of Quotients, Torsion Theories}

Through the paper, a ring is an associative ring with unit. By a
module we mean a right module unless otherwise specified. We adopt
the usual definitions of the injective envelope $E(M)$ of a module
$M$, the class of essential and dense submodules (e.g. definitions
3.31, 3.26, 8.2. \cite{Lam}), and the maximal right (left) ring of quotients
$\Qmax(R)$ ($\Qlmax(R)$) of a ring $R$ (sections 13B and 13C in \cite{Lam}).
If $\Qmax(R)=\Qlmax(R),$ we write $\Qlrmax(R)$ for $\Qmax(R)=\Qlmax(R).$

$\Qmax(R) \subseteq E(R)$ in general. If $R$ is right nonsingular,
the notions of dense and essential ideal are the same, $\Qmax(R)$
is equal to $E(R)$ and is von Neumann regular (Theorem 13.36 in
\cite{Lam}).

Let $S$ be a ring extension of $R.$ $S$ is a {\em general right
ring of quotients} if $R$ is dense in $S$ as a right $R$ module
(Definition 13.10 in \cite{Lam}). If $S$ is any general right
ring of quotients, then there is unique embedding of $S$ into
$\Qmax(R)$ that is identity on $R$ (Theorem 13.11, \cite{Lam}).

A {\em torsion theory} for $R$ is a pair $\tau = (\te, \f)$ of
classes of $R$-modules such that $\te$ and $\f$ are maximal classes having the property that
$\homo_R(T,F)=0,$ for all $T \in \te$ and $F \in \f.$
The modules in $\te$ are called {\em torsion modules} for $\tau$ and the modules in $\f$ are called
{\em torsion-free modules} for $\tau$.

A given class $\te$ is a torsion class of a torsion theory if an only if it is closed
under quotients, direct sums and extensions. A class $\f$ is a torsion-free class of a torsion theory if it is
closed under taking submodules, isomorphic images, direct products and extensions (see Proposition 1.1.9 in \cite{Bland}).

If $\tau_1 = (\te_1, \f_1)$ and $\tau_2 = (\te_2, \f_2)$ are two
torsion theories, we say that $\tau_1$ is {\em smaller} than
$\tau_2$ $(\tau_1\leq\tau_2$)  iff $\te_1\subseteq\te_2,$
equivalently $\f_1\supseteq\f_2.$

For every module $M$, the largest submodule of $M$ that belongs to $\te$ is called the
{\em torsion submodule} of $M$ and is denoted by $\te M$ (see Proposition 1.1.4 in \cite{Bland}). The
quotient $M/\te M$ is called the {\em torsion-free quotient} and
is denoted by $\f M.$ If $K$ is a
submodule of $M,$ the {\em closure} $\cl_{\tau}^M(K)$ of $K$ in $M$ with respect to the
torsion theory $\tau$ is largest submodule of $M$ such that $\cl_{\tau}^M(K)/K$ is torsion module (equivalently $M/\cl_{\tau}^M(K)$ is torsion-free).

A torsion theory $\tau = (\te, \f)$ is {\em hereditary} if the
class $\te$ is closed under taking submodules (equivalently torsion-free class is closed under
formation of injective envelopes, see Proposition 1.1.6,
\cite{Bland}). The largest torsion theory in which a given class of
injective modules is torsion-free (the torsion theory {\em cogenerated} by that class) is hereditary. Some authors
(e.g. \cite{Golan}, \cite{Lambek}) consider just hereditary torsion theories.
A torsion theory $\tau = (\te, \f)$ is {\em faithful}
if $R\in \f.$

The notion of Gabriel filter (terminology from \cite{Bland}) or
Gabriel topology (as is called in \cite{Stenstrom}) is equivalent
to the notion of hereditary torsion theory.

If $M$ is a $R$-module with submodule $N$ and $m$ an element of
$M,$ denote $\{r\in R\; | \; mr\in N\}$ by $(N : m).$
A {\em Gabriel filter (or Gabriel topology)} $\ef$ on a ring $R$
is a nonempty collection of right $R$-ideals such that
\begin{enumerate}
\item If $I\in \ef$ and $r\in R,$ then $(I:r)\in \ef.$

\item If $I\in \ef$ and $J$ is a right ideal with $(J:r)\in
\ef$ for all $r\in I,$ then $J\in \ef$.
\end{enumerate}

If $\tau$ is a hereditary torsion theory, the collection of right
ideals
$\{ I | R/I$ is a torsion module $\}$ is a Gabriel
filter $\ef_{\tau}.$ Conversely, if $\ef$ is a Gabriel
filter, then the class of modules $\{ M |  (0:m)$ is in
$\ef,$ for every $m\in M\}$ is a torsion class of a
hereditary torsion theory $\tau(\ef)$.The details can
be found in \cite{Bland} or \cite{Stenstrom}.

We recall some important examples of torsion theories.
\begin{exmp}
{\em

(1) The torsion theory cogenerated by the injective
envelope $E(R)$ of $R$ is called the {\em Lambek torsion theory}.
It is hereditary, as it is cogenerated by an injective module, and
faithful. Moreover, it is the largest hereditary faithful torsion
theory. The Gabriel filter of this torsion theory is the set of
all dense right ideals (see Proposition VI 5.5, p. 147 in
\cite{Stenstrom}).

(2) The class of nonsingular modules over a ring $R$ is closed
under submodules, extensions, products and injective envelopes.
Thus, it is a torsion-free class of a hereditary torsion theory.
This torsion theory is called the {\em Goldie torsion theory}. It
is larger than any hereditary faithful torsion theory (see Example
3, p. 26 in \cite{Bland}). So, the Lambek torsion theory is
smaller than the Goldie's. If $R$ is right nonsingular, the Lambek and Goldie torsion theories
coincide (see \cite{Bland} p. 26 or \cite{Stenstrom} p. 149).

(3) If $R$ is a right Ore ring with the set of regular elements $T$ (i.e.,
$rT \cap tR \neq 0,$ for every $t \in T$ and $r\in R$), we can
define a hereditary torsion theory by the condition that a right
$R$-module $M$ is a torsion module iff for every $m\in M$, there
is a nonzero $t\in T$ such that $mt =0.$ This torsion theory is
called the {\em classical torsion theory of a right Ore ring}. It
is hereditary and faithful.

(4) Let $R$ be a subring of a ring $S$. The collection of all $R$-modules $M$ such that $M\otimes_R S = 0$
is closed under quotients, extensions and direct
sums. Moreover, if $S$ is flat as a left $R$-module, then this
collection is closed under submodules and, hence, defines a
hereditary torsion theory. In this case we denote this torsion
theory by $\tau_S.$ From the definition of $\tau_S$ it follows that the torsion submodule of $M$ is the kernel of the
natural map $M\rightarrow M \otimes_R S$ and that
all flat modules are $\tau_S$-torsion-free. Thus, $\tau_S$ is faithful.
If $R$ is a right Ore ring, then $\tau_{\Qcl(R)}$ is the classical
torsion theory.
}
\label{Examples}
\end{exmp}

\subsection{Right Rings of Quotients}\label{subsection_on_right_rings_of_quotients}

If $\tau$ is a hereditary torsion theory with Gabriel filter $\ef
= \ef_{\tau}$ and $M$ is a right $R$-module, define:
\[M_{(\ef)} = \dirlim_{I\in\ef}\homo_R(I, M).\] In section 1 of
chapter 9 of \cite{Stenstrom} it is shown that $R_{(\ef)}$ has a
ring structure and that $M_{(\ef)}$ has a structure of a right
$R_{(\ef)}$-module.

Consider the map $\phi_M:M\rightarrow M_{(\ef)}$ obtained by
composing the isomorphism $M\cong\homo_R(R, M)$ with the map
$\homo_R(R, M)\rightarrow\dirlim\homo_R(I, M)$ given by $f\mapsto
f|_{I}.$ This $R$-homomorphism defines a left exact functor $\phi$ from the category of right $R$-modules to
the category of right $R_{(\ef)}$-modules.

\begin{lem}
\begin{enumerate}
\item $\te M = \ker (\phi_M: M \rightarrow M_{(\ef)}).$

\item $\te M = M$ if and only if $M_{(\ef)}=0.$

\item {\em coker}$\phi_M$ is a $\tau$-torsion module.
\end{enumerate}
\label{Lemma_M_ef}
\end{lem}

For details of the proof see Lemmas
IX 1.2, 1.3 and 1.5, p. 196 in \cite{Stenstrom}.

By parts 2. and 3. of Lemma \ref{Lemma_M_ef},
$(M/\te M)_{(\ef)}= (M_{(\ef)})_{(\ef)}.$ The {\em module
of quotients} $M_{\ef}$ of $M$ with respect to $\tau$ is defined as
\[M_{\ef} = (M_{(\ef)})_{(\ef)} = (M/\te M)_{(\ef)} =
\dirlim_{I\in \ef}\homo_R(I, M/\te M).\] The ring structure on
$R_{\ef}$ and the $R_{\ef}$-module structure on $M_{\ef}$ are
induced from corresponding structures on $R_{(\ef)}$ and
$M_{(\ef)}.$ The ring $R_{\ef}$ is called the {\em
right ring of quotients with respect to the torsion theory
$\tau.$} In \cite{Lambek}, there is an equivalent approach to the
notion of the module of quotients: $M_{\ef}$ is defined as closure
of $M/\te M$ in $E(M/\te M)$ with respect to $\tau.$ From this approach it readily follows that $M_{\ef}$ is
torsion-free as it is a submodule of an injective envelope of a
torsion-free module. Also, if
$\tau$ is faithful, then $R_{\ef}=\cl^{E(R)}_{\tau}(R).$

For every $M$, we have canonical homomorphism of $R$-modules $f_M:
M\rightarrow M_{\ef}.$ In particular, $f_R: R\rightarrow R_{\ef}$
is a ring homomorphism. The kernel of $f_M$ is he torsion module
$\te M$ for every module $M$ (see \cite{Stenstrom}, p. 197).

\begin{exmp}
{\em

(1) Since $\Qmax(R)=\dirlim\homo_R(I, R)$ where the limit is taken
over the family of dense ideals $I,$ $\Qmax(R)$ is the right ring
of quotients with respect to the Lambek torsion theory.

(2) Let $\ef_G$ be the filter of the Goldie torsion theory $\tau_G=(\te, \f).$
If $M$ is nonsingular,
its module of quotients $M_{\ef_G}$ is the injective envelope
$E(M)$ (see Propositions IX 2.5 and 2.7, Lemma IX 2.10 and Proposition IX
2.11 in \cite{Stenstrom}). For any $M$, $M_{\ef_G} =
\dirlim\homo_R(I, M)$ (Propositions IX 1.7 and
VI 7.3 in \cite{Stenstrom}), so $\dirlim\homo_R(I, M)= M_{\ef_G} =
\dirlim\homo_R(I, M/\te M) = (\f M)_{\ef_G}= E(\f M).$

If $R$ is right
nonsingular, $R_{\ef_G}=E(R)=\Qmax(R).$

(3) If $R$ is right Ore, the right ring of quotients with respect
to classical torsion theory (see part (3) of Example \ref{Examples})
is the classical right ring of quotients $\Qcl(R)$ (see Example 2,
ch. IX, p. 200 of \cite{Stenstrom}).
}
\end{exmp}

Let $S$ be a ring extension of $R.$ $S$ is a {\em right ring of
quotients} if $S=R_{\ef}$ for some Gabriel filter $\ef$ of a
hereditary torsion theory $\tau.$ In \cite{Lambek}, Lambek studies
the necessary and sufficient conditions for a ring extension $S$
to be a right ring of quotients.

If $\tau$ is hereditary and faithful with Gabriel filter $\ef$,
then $R_{\ef}$ can be embedded in $\Qmax(R)$ as $\tau$ is
contained in the Lambek torsion theory (see (1) of Example
\ref{Examples}). Since $R$ is dense in $\Qmax(R),$ then $R$
is dense in $R_{\ef}$ as well. So, a right ring of quotients
$R_{\ef}$ is also a general right ring of quotients if $\tau$ is faithful.

\subsection{Perfect Right Rings of Quotients}

Recall that the ring homomorphism $f:R\rightarrow S$ is
called a {\em ring epimorphism} if for all rings $T$ and
homomorphisms $g,h: S\rightarrow T,$ $gf = hf$ implies $g=h.$

\begin{prop} $f:R\rightarrow S$ is a ring epimorphism if and only if the canonical
map $S\otimes_R S\rightarrow S$ is bijective. \label{epimorphism}
\end{prop}

For proof see Proposition XI 1.2, p. 226 in \cite{Stenstrom}.

The situation when $S$ is flat as left $R$-module is of special
interest. There is a characterization of such epimorphisms due to
Popescu and Spircu (\cite{PopSpirc}).

\begin{thm} For a ring homomorphism $f:R\rightarrow S$ the following conditions are equivalent.
\begin{enumerate}
\item $f$ is a ring epimorphism and $S$ is flat as a left
$R$-module.

\item The family of right ideals $\ef=\{I | f(I)S=S\}$ is a
Gabriel filter, there is an isomorphism $g: S\cong R_{\ef}$ and
$g\circ f$ is the canonical map $R\rightarrow R_{\ef}.$

\end{enumerate}
\label{PerfectQuotient}
\end{thm}

The proof can also be found in \cite{Stenstrom}, p. 227.

If $f:R\rightarrow S$ satisfies the equivalent conditions of this
theorem, $S$ is called a {\em perfect right ring of quotients,} a {\em flat epimorphic extension} of $R,$
a {\em perfect right localization of $R$}
or a {\em flat epimorphic right ring of quotients of $R$.}

A hereditary torsion theory $\tau$ with Gabriel filter $\ef$ is
called {\em perfect} if the right ring of quotients $R_{\ef}$ is
perfect and $\ef=\{I| f_I(I)R_{\ef}=R_{\ef}\}$. The Gabriel filter
$\ef$ is called {\em perfect} in this case.

The perfect filters have a nice description. For a Gabriel filter
$\ef,$ let us look at the canonical maps $i_M: M\rightarrow
M\otimes_R R_{\ef}$ and $f_M: M\rightarrow M_{\ef}.$ There is a
unique $R_{\ef}$-map $F_M: M\otimes_R R_{\ef}\rightarrow
M_{\ef}$ given by $f_M = F_M i_M.$ The perfect filters are
characterized by the property that the map $F_M$ is an isomorphism
for every module $M.$ Moreover, the following holds.

\begin{thm} The following properties of a Gabriel filter $\ef$ are
equivalent.
\begin{enumerate}
\item $\ef$ is perfect.

\item The functor $q$ mapping the category of $R$-modules to the
category of $R_{\ef}$-modules given by $M\mapsto M_{\ef}$ is exact
and preserves direct sums.

\item $\ef$ has a basis consisting of finitely generated ideals
and the functor $q$ is exact.

\item The kernel of $i_M: M\rightarrow M\otimes_R R_{\ef}$ is a
torsion module in the torsion theory determined by $\ef$ for every
module $M.$

\item The map $F_M:M\otimes_R R_{\ef}\rightarrow M_{\ef}$ is an
isomorphism for every $M.$
\end{enumerate}
\label{perfect_filter}
\end{thm}

The proof can be found in \cite{Stenstrom} (Theorem XI 3.4, p.
231). Note that the functor $q$ from parts (2) and (3) is always
left exact.

This theorem establishes a one-to-one correspondence between the
set of perfect filters $\ef$ on $R$ and the perfect right rings of
quotients given by $\ef\mapsto R_{\ef}$ with the inverse $S\mapsto
\{I|f(I)S=S\}$ for $f:R\rightarrow S$ epimorphism that makes $S$ a
flat $R$-module.

From parts (4) and (5), it follows that if $\ef$ is a perfect
filter of torsion theory $\tau$, then $\tau$ is faithful because
then the torsion submodule of $R$ is isomorphic to $\tor^R_1(R,
R_{\ef}/R)$ which is 0 (see part (1) of Lemma \ref{Lemma_M_ef} and part (4) of Example \ref{Examples}). Thus,
if $S$ is a perfect right ring of quotients, then $R\subseteq
S\subseteq \Qmax(R).$

\subsection{The Total Right Ring of Quotients}

We further refine the introduced notions by considering the
maximal perfect right ring of quotients. Every ring
has a maximal perfect right ring of quotients, unique up to
isomorphism (Theorem XI 4.1, p. 233, \cite{Stenstrom}). It is
called {\em total right ring of quotients} (also maximal perfect
right localization of $R$, maximal flat epimorphic right ring of
quotients of $R$, right perfect hull, right flat-epimorphic hull).
We shall use the same notation as in \cite{Stenstrom} and denote it by $\Qtot(R).$
Other notations used in the literature include
epi$(R)$ and $M(R).$

In Theorem XI 4.1, p. 233, \cite{Stenstrom}, $\Qtot(R)$ is
obtained as the directed union of the family of all subrings of
$\Qmax(R)$ that are perfect right rings of quotients of $R.$ The
approaches in \cite{Findlay}, \cite{Knight}, \cite{Lazard}, and
\cite{PopSpirc} are all equivalent and involve the construction of
$\Qtot(R)$ as a direct limit. In \cite{Morita3}, Morita constructs
$\Qtot(R)$ differently than \cite{Findlay}, \cite{Knight},
\cite{Lazard} or \cite{PopSpirc}. If $M$ is a right $R$-module,
let us consider
\[\ef_t(M)=\{ I | I\mbox{ is a right ideal of }R\mbox{ and
}(I:r)M=M\mbox{ for all }r\in R\}.\] In Lemma 1.1 of
\cite{Morita3}, Morita shows that this is a Gabriel filter of a
hereditary torsion theory.

In Theorem 3.1 of \cite{Morita3}, Morita shows that a ring
homomorphism $f:R\rightarrow S$ is a ring epimorphism with $S$
flat as a left $R$-module if and only if $S$ is the right ring of
quotients of $R$ with respect to the Gabriel filter $\ef_t(S).$
In this case $S=\{s\in S | (R:sr)S=S$ for every $r\in R\}.$

Motivated by this result Morita considers the set
\[S'=\{s\in S | (R:sr)S=S\mbox{ for every }r\in R\}\]
for a ring extension $S$ of $R.$ By Theorem 3.1 of \cite{Morita3}, $S$ is flat
epimorphic extension if and only if $S=S'.$ In Lemma 3.2 of
\cite{Morita3}, Morita proves that $S'$ is a subring of $S$ that
contains $R$ for a ring extension $S$ of $R.$ In Corollary 3.4 of
\cite{Morita3}, he shows that there exist the largest flat
epimorphic extension of $R$ that is contained in a given extension
$S.$ After proving this corollary, Morita also sketches the idea
of the alternative proof (passage following the proof). We are
interested in this alternative proof. The outline of the proof is
the following.

Let $S^{(0)}=S.$ If $\alpha$ is a successor ordinal
$\alpha=\beta+1,$ then $S^{(\alpha)}=(S^{(\beta)})'.$ If $\alpha$
is a limit ordinal, let $S^{(\alpha)}=\bigcap_{\beta<\alpha}
S^{(\beta)}.$ Morita claims that there is an ordinal $\gamma$ such
that $S^{(\gamma)}=(S^{(\gamma)})'=S^{(\gamma+1)}.$ This is true because if $S^{(\gamma+1)}$ is strictly contained in $S^{(\gamma)}$ for every ordinal $\gamma,$ then $|S|\geq |S-S^{(\gamma)}|\geq |\gamma|$ for every ordinal $\gamma$ which is a contradiction. If
$S^{(\gamma)}=S^{(\gamma+1)},$ then $S^{(\gamma)}$ is flat
epimorphic extension of $R$ by Theorem 3.1 in \cite{Morita3}. To
see that $S^{(\gamma)}$ is the largest flat epimorphic extension
contained in $S,$ take $T$ to be any flat epimorphic extension
such that $T\leq S.$ Then $T'=T\leq S'$ so it is easy to see that
$T$ is contained in all extensions $S^{(\alpha)}$ for every
ordinal $\alpha.$ Hence, $T\leq S^{(\gamma)}.$

$S=\Qmax(R)$ is the case of special interest. In this case,
this construction gives us $\Qtot(R)$ (see last paragraph of
Section 3 in \cite{Morita3}). In the rest of the paper, we shall refer to this construction of
$\Qtot(R)$ as Morita's construction.

\begin{exmp}
{\em
(1) If $R$ is regular, then $R=\Qtot(R)$ by Example 1 and
Proposition XI 1.4, p. 226 in \cite{Stenstrom}.

(2) If $R$ is right Ore, then $\Qcl(R)\subseteq\Qtot(R).$ If
$\Qcl(R)$ is regular, then $\Qcl(R)=\Qtot(R)$ (Example 2, ch. XI,
p. 235, \cite{Stenstrom}).

(3) If $R$ is right noetherian and right hereditary (in particular
if $R$ is semisimple), then $\Qmax(R)=\Qtot(R)$ (Example 3, ch.
XI, p. 235, \cite{Stenstrom}) If $R$ is also commutative, then
$\Qlrcl(R)=\Qlrmax(R)=\Qlrtot(R).$ \label{Examples_of_perfect}

}
\end{exmp}

\section{Construction of $\Qtot(R)$ for a class of rings}
\label{Section_construction}

In this section, we consider a class of rings for which the
Gabriel filter from Morita's construction at step $\alpha$ is
exactly the Gabriel filter of the torsion theory obtained by
tensoring with $\Qmax(R)^{(\alpha)}$ (see part (4) of Example
\ref{Examples}) for all ordinals $\alpha.$ First, we need
the following lemma.

\begin{lem} Let $\tau=(\te, \f)$ be a hereditary torsion theory with Gabriel filter
$\ef$ such that its right ring of quotients $R_{\ef}$ is flat as
left $R$-module.
\begin{itemize}
\item[1.] The torsion theory $\tau_{R_{\ef}}$ (introduced in (4) of Example \ref{Examples}) is smaller than $\tau.$ If $\tau$ is
faithful, the right ring of quotients of $\tau_{R_{\ef}}$ is
contained in $R_{\ef}.$

\item[2.] $\tau = \tau_{R_{\ef}}$ if and only if $\tau$ is
perfect.

\item[3.] If $R_{\ef}$ is a perfect right ring of quotients then
the torsion theory $\tau_{ R_{\ef}}$ is perfect.
\end{itemize}
\label{t_is_in_T}
\end{lem}

Note that in the last part of this lemma, it is possible to have
$R_{\ef}$ (and $\tau_{R_{\ef}}$) perfect without $\tau$ being
perfect. We illustrate this situation in Example
\ref{Example_classC}.

\begin{pf} 1. Denote $\tau_{R_{\ef}}$
with $(t, p).$ We will show that $t\subseteq \te.$ Let $M$ be any
right $R$-module. $t M$ is the kernel of $i_M: M\rightarrow
M\otimes_R R_{\ef}$ (see part (4) of Example \ref{Examples}). It is
contained in $\ker (f_M: M\rightarrow M_{\ef}).$ But $\ker f_M$ is
$\te M.$ Thus, $t M\subseteq \te M.$

Let $S$ be the right ring of right quotients of torsion theory
$(t,p).$ $(t,p)$ is faithful so $S=\dirlim\homo_R( I, R)$ where
the limit is taken over the right ideals $I$ that are in the
Gabriel filter of $(t,p).$ Since $\tau$ is faithful as well,
$R_{\ef}=\dirlim\homo_R(I,R),$ $I\in \ef.$ But the filter
corresponding to $(t,p)$ is contained in $\ef$ and so $S\subseteq
R_{\ef}.$

2. If $t M=\te M,$ then condition (4) from Theorem
\ref{perfect_filter} holds so $\tau$ is perfect.
Conversely, if $\tau$ is perfect and $M$ is a torsion with respect
to $\tau,$ then $M_{\ef}=0$ by part (2) of Lemma \ref{Lemma_M_ef}.
But $F_M$ is an isomorphism by condition (5) of Theorem
\ref{perfect_filter}, so $M\otimes_R R_{\ef}=0.$ Hence, $M$ is
torsion in $(t, p)$ by part (4) of Example \ref{Examples} so the two
torsion theories coincide.

3. If $R_{\ef}$ is perfect, then it is a right ring of quotients
of a perfect torsion theory (not necessarily $\tau$). That torsion
theory is equal to $\tau_{R_\ef}$ by part 2. So, $\tau_{R_{\ef}}$
is perfect.
\end{pf}

The idea of our construction is to start by checking if Lambek
torsion theory is perfect. Denote its right ring of quotients
$\Qmax(R)$ by $Q_0.$ If it is perfect, $Q_0=\Qtot(R).$ If not, we
consider the strictly smaller torsion theory $\tau_{Q_0}.$ If it
is perfect, its right ring of quotients $Q_1$ is $\Qtot(R).$ If
not, we consider the strictly smaller torsion theory $\tau_{Q_1}$
and continue inductively. If the construction does not end after
finitely many steps, we consider $Q_{\omega}$ to be show the
intersection of the rings $Q_n,$ $n\geq 0,$ and proceed
inductively.

The only thing we need to insure in order to be able to define the
above torsion theories and their rings of quotients is that the
defined ring extensions of $R$ are flat as left $R$-modules. Thus,
we impose the following condition on $R:$

\begin{itemize}
\item[(C)] Every subring of $\Qmax(R)$ that contain $R$ is flat as
a left $R$-module.
\end{itemize}
Under this condition, let us prove that the above described idea
works.

{\bf Step 0.} Denote the Lambek torsion theory by $\tau_0$, its
filter, the set of all dense right ideals by $\ef_0,$ and its
right ring of quotients, $\Qmax(R)$ by $Q_0.$

Check if $\tau_0$ is perfect. Note that, if $R$ is right
nonsingular, this is equivalent to the condition that $\Qmax(R)$
is semisimple by Proposition XI 5.2 and Example 2, p. 237 in
\cite{Stenstrom}. If $\tau_0$ is perfect, then
$\Qtot(R)=Q_0=\Qmax(R)$ by (3) of Examples
\ref{Examples_of_perfect} and the construction is over. If not, go
to next step.

{\bf Inductive step.} Let us suppose that we constructed the
torsion theory $\tau_{\alpha}$ with Gabriel filter $\ef_{\alpha}$
and the right ring of quotients $Q_{\alpha}.$ Then, we define
\[\tau_{\alpha+1}=\tau_{Q_{\alpha}},\;\;\;\ef_{\alpha+1}=\mbox{ Gabriel filter
corresponding to
}\tau_{\alpha+1},\;\;\;Q_{\alpha+1}=R_{\ef_{\alpha+1}}.\] Here we
are using condition (C) in order for $\tau_\alpha$ to be
hereditary.

If $\alpha$ is a limit ordinal and the rings $Q_{\beta}$ for
$\beta< \alpha$ are constructed, then define
\[\tau_{\alpha}=\bigcap_{\beta<\alpha}\tau_{\beta},\;\;\;\ef_{\alpha}=\mbox{ Gabriel filter
corresponding to }\tau_{\alpha}=\bigcap_{\beta<\alpha}
\ef_{\beta},\;\;\;Q_{\alpha}=R_{\ef_{\alpha}}.\]

Note that in this case
$Q_{\alpha}=\bigcap_{\beta<\alpha}Q_{\beta}.$ One direction
follows since
$\ef_{\alpha}\subseteq\bigcap_{\beta<\alpha}\ef_{\beta}.$ To prove
the other direction, let us note that
$Q_{\beta}=\cl_{\tau_{\beta}}^{E(R)}(R)$ as every $\tau_{\beta}$
is faithful. Then $(\bigcap Q_{\beta})/R$ is torsion in
$\tau_{\beta}$ for every $\beta<\alpha$ as it is a submodule of
torsion module $Q_{\beta}/R=\cl_{\tau_{\beta}}^{E(R)}(R)/R.$  So,
$\bigcap Q_{\beta}$ has to be contained in the closure
$\cl_{\tau_{\alpha}}^{E(R)}(R)=Q_{\alpha}.$

Let us note also that $Q_{\alpha}/R$ is a torsion module in
$\tau_{\alpha}$ as is the cokernel of map $R\hookrightarrow
Q_{\alpha}$ (see part (3) of Lemma \ref{Lemma_M_ef}).

\begin{lem} Let $\beta<\alpha.$
\begin{enumerate}

\item $\tau_{\alpha}\subseteq \tau_{\beta}$ and
$Q_{\alpha}\subseteq Q_{\beta}.$

\item  $Q_{\beta}/Q_{\alpha}$ is torsion module in $\tau_{\beta}$
and torsion-free module in $\tau_{\alpha}.$

\item $Q_{\alpha}\otimes_R Q_{\beta}=R\otimes_R Q_{\beta}\cong
Q_{\beta}.$

\item $\Qtot(R)\subseteq Q_{\alpha}.$

\item $\tau_{\beta}=\tau_{\alpha}$ if and only if $\tau_{\beta}$
is perfect.

\item $Q_{\alpha}$ is perfect right ring of quotients if and only
if $Q_{\alpha}= \Qtot(R).$

\item If $\tau_{\alpha}$ is perfect, then $Q_{\alpha}$ is perfect.
If $Q_{\alpha}$ is perfect, then $\tau_{\alpha+1}$ is perfect.
\end{enumerate}
\label{Induction_Works}
\end{lem}

\begin{pf}

(1) This is part 1. of Lemma \ref{t_is_in_T} for $\alpha$
successor ordinal and definition of $\tau_{\alpha}$ for $\alpha$
limit ordinal.

(2) $Q_{\beta}/Q_{\alpha}$ is a quotient of $Q_{\beta}/R.$
$Q_{\beta}/R$ is torsion in $\tau_{\beta}$ and then so is
$Q_{\beta}/Q_{\alpha}.$

$Q_{\beta}/Q_{\alpha}$ is a submodule of $E(R)/Q_{\alpha}.$ But
$Q_{\alpha}=\cl_{\tau_{\alpha}}^{E(R)}(R)$ so $E(R)/Q_{\alpha}$ is
torsion-free in $\tau_{\alpha}.$ Hence, the submodule
$Q_{\beta}/Q_{\alpha}$ is torsion-free in $\tau_{\alpha}$ as well.

(3) $\beta<\alpha$ implies $\beta+1\leq \alpha$. $Q_{\alpha}/R\leq
Q_{\beta+1}/R$ is torsion in $\tau_{\beta+1}.$ Thus,
$Q_{\alpha}/R\otimes_R Q_{\beta}=0.$ Since $Q_{\beta}$ is flat, we
have that $Q_{\alpha}\otimes_R Q_{\beta}=R\otimes_R Q_{\beta}\cong
Q_{\beta}.$

(4) We show this by induction on $\alpha.$ If $\alpha=0,$
$\Qtot(R)\subseteq\Qmax(R)=Q_0$ as $\Qtot(R)$ is a general right
ring of quotients. Suppose that it holds for all ordinals less
than $\alpha.$ If $\alpha$ is a limit ordinal, the claim easily
follows. Let $\alpha$ be a successor ordinal of $\beta.$ Let $q\in
\Qtot(R).$ Then $q$ can be represented as a map $I\rightarrow R$
for some right ideal $I$ with $I\Qtot(R)=\Qtot(R)$ by part (2) of
Theorem \ref{PerfectQuotient}. So, $1=\sum r_i q_i$ for some
$r_i\in I$ and $q_i\in \Qtot(R),$ $i=1,\ldots,m$ for some $m.$ By
induction hypothesis, $q_i$ is in $Q_{\beta}.$ Thus
$Q_{\beta}\subseteq IQ_{\beta}$ and so $IQ_{\beta}=Q_{\beta}.$
Hence, $q$ is in the right ring of quotients with respect to
$\tau_{Q_{\beta}}$ which is $Q_{\alpha}.$

(5) Since $\beta<\alpha$ implies $\beta+1\leq \alpha$,
$\tau_{\beta}=\tau_{\alpha}$ implies
$\tau_{\beta}=\tau_{\beta+1}.$ Then $\tau_{\beta}$ is perfect by
part 2. of Lemma \ref{t_is_in_T}. Conversely, if $\tau_{\beta}$ is perfect, then
$\tau_{\beta}=\tau_{\beta+1}$ (again by part 2. of Lemma
\ref{t_is_in_T}) so $\tau_{\beta}=\tau_{\alpha}$ for all
$\alpha>\beta.$

(6) If $Q_{\alpha}$ is perfect, $Q_{\alpha}$ is contained it
$\Qtot(R)$ by definition of $\Qtot(R).$ Since the converse always
holds by part (4), we have that $Q_{\alpha}=\Qtot(R).$ The
converse is clear.

(7) The first part follows from Theorem \ref{perfect_filter} and the
second part from part 3. of Lemma \ref{t_is_in_T}.
\end{pf}

From part (7), we see that $\tau_{\alpha}$ being perfect implies
that $Q_{\alpha}$ is perfect as well. The converse does not hold
(see Example \ref{Example_classC}). Also, if $Q_{\alpha}$ is
perfect, $\tau_{\alpha+1}$ is perfect as well but the converse
does not have to hold (see Example \ref{Q0_not_Qtot}).

{\bf Getting $\Qtot(R).$} Ordinal $\alpha$ such that
$Q_{\alpha}=Q_{\alpha+1}$ has to exist by the same argument as the one used in the proof of Morita's
construction. If $Q_{\alpha}=Q_{\alpha+1},$ then $Q_{\alpha}\otimes_R Q_{\alpha}=Q_{\alpha+1}\otimes_R Q_{\alpha}
\cong Q_{\alpha}$ by part (3) of Lemma \ref{Induction_Works}. Thus $Q_{\alpha}$ is perfect by Proposition \ref{epimorphism}. Then $Q_{\alpha}=\Qtot(R)$ by part (6) of Lemma \ref{Induction_Works}.

The next proposition shows that Morita's construction coincides
with our construction if the ring $R$ satisfies condition (C).

\begin{prop} If $R$ is a ring that satisfies (C), then for
$Q=\Qmax(R),$ \[Q_{\alpha}=Q^{(\alpha)} \mbox{ for all }\alpha.\]
\label{Morita=mine_whenC}
\end{prop}

\begin{pf} $Q_0=Q^{(0)}$ as both
are $\Qmax(R).$ Let us proceed by induction. Assume that
$Q_{\alpha}=Q^{(\alpha)}.$ Recall that $Q_{\alpha+1}$ is the right
ring of quotients with respect to the Gabriel filter
$\ef_{\alpha+1}=\{I | IQ_{\alpha}=Q_{\alpha}\}.$ $Q^{(\alpha+1)}$
is the right ring of quotients with respect to the Gabriel filter
$\ef_t(Q^{(\alpha)})=\{I | (I:r)Q^{(\alpha)}=Q^{(\alpha)}$ for all
$r\in R \}$ by Theorem 4.1 of \cite{Morita3}. Clearly if $I$ is a
right ideal in $\ef_t(Q^{(\alpha)}),$ then $(I:
1)Q^{(\alpha)}=Q^{(\alpha)}$ and so $IQ_{\alpha}=Q_{\alpha}.$
Conversely, if $I$ is in $\ef_{\alpha+1},$ then $(I:r)$ is in
$\ef_{\alpha+1}$ for any $r\in R$ by property (1) of Gabriel
filter (see the definition of Gabriel filter in Section
\ref{Section_general_quotients}). Since we assume that
$Q_{\alpha}=Q^{(\alpha)},$ then $I\in\ef_t(Q^{(\alpha)}).$

If $\alpha$ is a limit ordinal and we assume that
$Q_{\beta}=Q^{(\beta)}$ for all $\beta<\alpha,$ then
$Q_{\alpha}=\bigcap Q_{\beta}=\bigcap Q^{(\beta)}=Q^{(\alpha)}.$
\end{pf}

\section{$\Qtot(R)$ of a Right
Semihereditary Ring $R$} \label{Section_semihereditary}

In this section, we consider the class of right semihereditary
rings to illustrate the benefits of using our construction when it
is possible to do so. Let us first prove the following lemma.

\begin{lem} For any $R$ that satisfies (C), the Gabriel filter $\ef_{\alpha}$ has a basis consisting of finitely generated right ideals
for every successor ordinal $\alpha$. \label{finite_basis_Lemma}
\end{lem}

\begin{pf} The statement of the lemma means that for every right ideal
$I$ in $\ef_{\alpha},$ there is finitely generated right ideal $J$
in $\ef_{\alpha}$ such that $J\subseteq I.$

Let $I\in \ef_{\alpha}.$ Since $\alpha$ is successor,
$\alpha=\beta+1$ for some $\beta.$ By construction, this means
that $I Q_{\beta}=Q_{\beta}.$ Then, there is $m$ and $r_i\in I,$
$q_i\in Q_{\beta},$ $i=1,\ldots, m$ such that $\sum r_i q_i=1.$

Let $J$ be the right ideal generated by $\{r_1, \ldots, r_m\}.$
Clearly, $J\subseteq I.$ $1=\sum r_i q_i\in JQ_{\beta}$ and so
$Q_{\beta}=J Q_{\beta}.$ Thus, $J$ is in $\ef_{\alpha}.$
\end{pf}

This lemma is the essential reason why it is better to consider
Gabriel filters $\ef_{\alpha}$ instead of $\ef_t(Q^{(\alpha)})$
when possible. In general, there is no reason for the filter
$\ef_t(Q^{(\alpha)})$ to have a basis consisting of finitely
generated ideals and the usefulness of the property is evident in
part (3) of Theorem \ref{perfect_filter}. On the other hand,
filters $\ef_{\alpha}$ do have this property for $\alpha$
successor by Lemma \ref{finite_basis_Lemma}. This property of
filters $\ef_{\alpha}$ will be essential when considering the
class of right semihereditary rings in the next theorem.

\begin{thm} If $R$ is right semihereditary, then $R$ satisfies (C) and \[\Qtot(R)=
Q_1.\] \label{Qtot_for_semihereditary}
\end{thm}

\begin{pf}
$\Qmax(R)$ is left flat for every right nonsingular and right
coherent ring $R$: a right coherent ring has a left flat right
ring of quotients with respect to the Goldie torsion theory
(Example 1, ch. XI, p. 233 \cite{Stenstrom}), and a right
nonsingular ring has equal Lambek and Goldie torsion theories, so
the Goldie right ring of quotients is the same as $\Qmax(R).$ (C)
is true if $R$ is, in addition, subflat. A ring is subflat if
every submodule of a left (equivalently right) flat $R$-module is
flat. Equivalently, all left (right) ideals are flat. Right
nonsingular, right coherent rings that are subflat are right
semihereditary (Theorem 2.10 in \cite{Sandomiersky} and Example 1,
p. 233 \cite{Stenstrom}). Converse also holds, if $R$ is right
semihereditary, then it is right nonsingular, right coherent and
subflat.

For the construction to end after the first step, it is sufficient
to show that the filter $\ef_1$ is perfect. We show that the
condition (3) from Theorem \ref{perfect_filter} is satisfied for
$\ef_1.$ By above lemma, $\ef_1$ has a basis of finitely generated
right ideals. But $R$ is right semihereditary so those ideals are
projective. Then the functor $q$ from condition (3) of Theorem
\ref{perfect_filter} is exact since any Gabriel filter $\ef$ with
basis consisting of projective right ideals has exact functor $q$
(Proposition XI 3.3, p. 230, \cite{Stenstrom}). So,
$Q_1=\Qtot(R).$
\end{pf}

This theorem provides us with a simple hands-on description of the
total right ring of quotients for $R$ right semihereditary:
\[\Qtot(R)=\{\; q\in \Qmax(R)\; |\; (R:q)\Qmax(R)=\Qmax(R)\;\}.\]
Let us consider the following examples of semihereditary rings.

\subsection{Example of a semihereditary ring with $Q_0=\Qtot(R),$ $\tau_0$ not perfect}
\label{Example_classC}

The class $\ce$ considered in \cite{Be2}, \cite{Lia2} and
\cite{Lia3} consists of certain finite Baer *-rings that are all
semihereditary (see Corollary 5 in \cite{Lia2}). All finite
$AW^*$-algebras (in particular all finite von Neumann algebras)
are in $\ce.$

A ring $R$ from $\ce$ has (left and right) maximal and classical
ring of quotients equal by Proposition 3 in \cite{Lia2} (let us
denote it by $Q$) and thus $\Qlrtot(R)$ is equal to $Q$ as well.
Thus, for this class of rings $Q_0=\Qtot(R).$ However, not all
rings in $\ce$ have $\tau_0$ perfect. In fact, part 3 of Theorem
23 in \cite{Lia2} says that $\tau_0=\tau_1$ (in notation used in
this paper) if and only if $Q$ is semisimple. This is equivalent
to the condition that $\tau_0$ is perfect by part 2 of Lemma
\ref{t_is_in_T}. The inequality $\tau_1\leq\tau_0$ can be strict
by Example 8.34 in \cite{Lu_book}. Note also that this is an
example of a ring with $\tau_0$ and $\tau_1$ different but with
the same right ring of quotients $Q_0.$ So, it is possible to have
the perfect $Q_0$ but not perfect $\tau_0.$

\subsection{Example of a semihereditary ring with $Q_0\neq Q_1= \Qtot(R)$}
\label{Q0_not_Qtot} Let $R=\{(a_n)\in
\Qset\times\Qset\times\ldots\;|\; (a_n)\mbox{ is eventually
constant }\}.$ $R$ is commutative so the left and right ring of
quotients coincide. $R$ is regular, so $\Qlrtot(R)=R.$
$\Qlrmax(R)=\Qset\times\Qset\times\ldots$ (Exercise 23, p. 328,
\cite{Lam}). As regular rings are semihereditary,
$Q_1=\Qlrtot(R)=R.$

This example also provides the evidence of a ring with $\tau_1$
perfect without $Q_0$ being perfect and a maximal ring of
quotients that is flat but not perfect.

Another example of a commutative ring with $\tau_0$ not perfect
can be found on page 332 in \cite{Schelter_Roberts}.

\subsection{Semihereditary Rings with $\Qmax(R)=\Qtot(R)$}
Let us mention some results related to the condition that
$\Qmax(R)=\Qtot(R).$ In general, this condition is weaker than the
condition that $\tau_0$ is perfect as we have seen in Example
\ref{Example_classC}.

In \cite{Goodearl}, Goodearl showed that for a right nonsingular
ring $R$, the following are equivalent:
\begin{itemize}
\item[i)] Every finitely generated nonsingular module can be
embedded in a free module.

\item[ii)] $\Qmax(R)=\Qltot(R).$
\end{itemize}
This result implies that the following two conditions on a right
nonsingular ring $R$ are equivalent:
\begin{enumerate}
\item Every finitely generated nonsingular module is projective.

\item $R$ is right semihereditary and $\Qmax(R)=\Qltot(R).$
\end{enumerate}
Also, if these conditions are satisfied then $R$ is also left
semihereditary and $\Qmax(R)=\Qtot(R).$ This result was first
shown in \cite{Cateforis1}.

In \cite{Evans_all}, Evans shows that the following conditions are
equivalent
\begin{itemize}
\item[(3)] $R$ is right semihereditary ring and
$\Qmax(R)=\Qtot(R)=\Qltot(R).$

\item[(4)] The matrix ring $M_n(R)$ is strongly Baer (every right
complement ideal is generated by an idempotent) for all $n.$
\end{itemize}
Evans calls the rings satisfying these equivalent conditions the
{\em right strongly extended semihereditary.} The rings from
Example \ref{Example_classC} are (left and right) strongly
extended semihereditary. The ring from Example \ref{Q0_not_Qtot}
is an example of a (left and right) semihereditary ring that is
not strongly extended semihereditary.

In \cite{Finkel}, Finkel Jones considers the notion of
$f$-projectivity. A module $M$ is said to be {\em $f$-projective}
if the inclusion of a finitely generated submodule of $M$ in $M$
factors through a free module. $f$-projectivity lies properly
between projectivity and flatness. Every finitely generated
$f$-projective module is projective. If $R_{\ef}$ is perfect ring
of quotients, then $R_{\ef}$ is $f$-projective by Proposition 2.1,
p. 1608 in \cite{Finkel}. Conversely, if $R_{\ef}$ is
a ring of quotients with respect to a faithful hereditary torsion
theory such that $R_{\ef}$ is $f$-projective, then $R_{\ef}$ is
perfect. Thus, the notion of $f$-projectivity also characterizes
the perfect right rings of quotients.

In \cite{Evans_all}, Evans uses the notion of $f$-projectivity to
further describe a class of right strongly extended semihereditary
rings. He proves that the following conditions are equivalent to
(3) and (4) above:
\begin{itemize}
\item[(5)] The class of $f$-projective modules is a torsion-free class of a
hereditary torsion theory.

\item[(6)] A module is $f$-projective if and only if it is
nonsingular.
\end{itemize}

\section{A Class of Rings for Which the Construction Ends After Countably Many Steps}
\label{Section_construction_stops}

Let $\omega$ denote the first infinite ordinal as usual.

\begin{prop} If $R$ satisfies condition (C) and
\begin{itemize}
\item[(C')] Every subring of $\Qmax(R)$ that contain $R$ is flat
as a right $R$-module,
\end{itemize}
then \[Q_{\omega}=\Qtot(R).\]

In particular, a commutative ring $R$ that satisfies condition (C)
has $Q_{\omega}=\Qtot(R).$ \label{C_and_C'}
\end{prop}

\begin{pf} Since $R$ satisfies (C), we know that $Q_{\omega}$ is flat
as a left $R$-module. Thus, to prove that it is perfect it is
sufficient to show that the canonical map $Q_{\omega}\otimes_R
Q_{\omega}\rightarrow Q_{\omega}$ is an isomorphism (by
Proposition \ref{epimorphism}). $Q_{\omega}\otimes_R
Q_{\omega}\leq Q_{\omega}\otimes_R Q_{n}$ as $Q_{\omega}$ is flat
as a right $R$-module by (C').

\[\begin{array}{rcll}
Q_{\omega}\otimes_R Q_{\omega} & \leq & \bigcap
(Q_{\omega}\otimes_R Q_n) & (\mbox{by what we showed above})\\
 & = & \bigcap(R\otimes_R Q_n)  & (\mbox{by part (3) of Lemma \ref{Induction_Works}})\\
 & = & R\otimes_R\bigcap Q_n & (\mbox{inverse limit commutes with }R\otimes_R\underline{\hskip0.3cm})\\
 & = & R\otimes_R Q_{\omega} & (\mbox{by definition of }Q_{\omega})\\
 & \cong & Q_{\omega} &
\end{array}\]

If $R$ is commutative, then $\Qmax(R)$ is commutative as well (see
Proposition 13.34 in \cite{Lam}). Thus condition (C) implies
condition (C') so the claim follows.
\end{pf}

Note that in the proof we really used much weaker assumption than
(C'). Namely, we just used that $Q_{\omega}$ is flat as right
$R$-module, not that every subring of $\Qmax(R)$ that contains $R$
is flat as right module. Thus, we obtain the following corollary.

\begin{cor} If $R$ is a ring that satisfies (C) and such that $Q_{\alpha}$ is flat as a right
$R$-module for some limit ordinal $\alpha,$ then
$\Qtot(R)=Q_{\alpha}.$
\end{cor}
To prove this, just replace $\omega$ with $\alpha$ and $n$ with
any $\beta<\alpha$ in the proof of Proposition \ref{C_and_C'}.

\section{Questions} We conclude by listing some interesting
questions and problems.

\begin{enumerate}
\item In \cite{Stenstrom}, p. 235, Stenstr\"om is asking for
necessary and sufficient conditions for $\Qmax(R)$ and $\Qtot(R)$
to be equal. Note that this is weaker than the condition for the
Lambek torsion theory to be perfect. The necessary and sufficient
condition for the Lambek torsion theory to be perfect is known:
$\tau_0$ is perfect if and only if $\Qmax(R)$ has no proper dense right
ideals (Proposition XI 5.2, p. 236, \cite{Stenstrom}). A ring $R$
satisfying this condition is called right {\em Kasch}. If $R$ is
hereditary and noetherian (Example 3, p. 235, \cite{Stenstrom}) or
commutative and noetherian (Example 4, p. 237, \cite{Stenstrom})
or nonsingular with finite uniform dimension (Gabriel's Theorem,
see Theorem 13.40 in \cite{Lam} or Theorem XII 2.5 in
\cite{Stenstrom}), $\Qmax(R)$ is known to be Kasch.

\item For any $n,$ find example of a ring $R$ such that
$Q_n=\Qtot(R)\neq Q_i$ for $i<n.$ Describe the rings satisfying
this condition.

\item Find example of a ring $R$ such that
$Q_{\omega}=\Qtot(R)\neq Q_n$ for all $n.$ Describe the rings
satisfying this condition.

\item In Example 4, p. 253 of \cite{Stenstrom}, Stenstr\"om is
asking how the type of Baer ring changes when taking the maximal
ring of quotients. With that in mind, it would also be natural to
ask how the type of Baer ring changes when taking the total ring
of quotients.
\end{enumerate}

\end{document}